\documentclass[12pt]{amsart}
\usepackage{amssymb,amsmath,amsthm}
\usepackage{fancyhdr}
\newtheorem{thm}{Theorem}[section]

\begin{document}

\title{Sobolev regularity of the Bergman and Szeg\"{o} projections in terms of $\overline{\partial}\oplus\overline{\partial}^{*}$ and $\overline{\partial}_{b}\oplus\overline{\partial}_{b}^{*}$}
\author{Emil J. Straube}
\thanks{Research supported in part by NSF grant DMS--2247175}

\address{Department of Mathematics, Texas A\&M University, College Station, Texas, USA}
\email{e-straube@tamu.edu}
\date{October 13, 2024}

\begin{abstract}
Let $\Omega$ be a smooth bounded pseudoconvex domain in $\mathbb{C}^{n}$. It is shown that for $0\leq q\leq n$, $s\geq 0$, the embedding $j_{q}: dom(\overline{\partial})\cap dom(\overline{\partial}^{*}) \hookrightarrow L^{2}_{(0,q)}(\Omega)$ is continuous in $W^{s}(\Omega)$--norms if and only if the Bergman projection $P_{q}$ is (see below for the modification needed for $j_{0}$). The analogous result for the operators on the boundary is also proved (for $n\geq 3$). In particular, $j_{1}$ is always regular in Sobolev norms in $\mathbb{C}^{2}$, notwithstanding the fact that $N_{1}$ need not be.
\end{abstract}

\maketitle

\section{Introduction and Results}\label{intro}

Let $\Omega$ be a smooth bounded pseudoconvex domain in $\mathbb{C}^{n}$. For $0\leq q\leq n$, denote by $L^{2}_{(0,q)}(\Omega)$ the usual Hilbert space of $(0,q)$ forms with coefficients in $L^{2}(\Omega)$. Analogously, for $0\leq q\leq (n-1)$, $L^{2}_{(0,q)}(b\Omega)$ denotes the Hilbert space of $(0,q)$ forms with square integrable coefficients on the boundary. See \cite{Boggess91, BoasStraube91b, BiardStraube17, ChenShaw01, Straube10} and their references for details on the $L^{2}$--setup of the $\overline{\partial}$ and $\overline{\partial}_{b}$ complexes. For $1\leq q\leq (n-1)$, denote by $j_{q}$ the embedding $j_{q}: dom(\overline{\partial})\cap dom(\overline{\partial}^{*}) \hookrightarrow L^{2}_{(0,q)}(\Omega)$, and by $P_{q}$ the (orthogonal) Bergman projection  $P_{q}: L^{2}_{(0,q)}(\Omega) \rightarrow \ker(\overline{\partial})\subseteq L^{2}_{(0,q)}(\Omega)$. When $q=0$,  $dom(\overline{\partial})\cap dom(\overline{\partial}^{*})$ is replaced by $dom(\overline{\partial})\cap\ker(\overline{\partial})^{\perp}$, and when $q=n$ by $dom(\overline{\partial}^{*})\cap\ker(\overline{\partial}^{*})^{\perp} = dom(\overline{\partial}^{*})$. Similarly, for $1\leq q\leq (n-2)$ denote by $j_{b,q}$ the embedding $j_{b,q}: dom(\overline{\partial}_{b})\cap dom (\overline{\partial}_{b}^{*}) \hookrightarrow L^{2}_{(0,q)}(b\Omega)$; for $q=0$, $dom(\overline{\partial}_{b})\cap dom (\overline{\partial}_{b}^{*})$ is replaced by $dom(\overline{\partial}_{b})\cap\ker(\overline{\partial}_{b})^{\perp}$, and for $q=(n-1)$ by $dom(\overline{\partial}_{b}^{*})\cap\ker(\overline{\partial}_{b}^{*})^{\perp}$. For $0\leq q\leq (n-2)$, the Szeg\"{o} projection $S_{q}$ is the orthogonal projection $S_{q}: L^{2}_{(0,q)}(b\Omega)\rightarrow \ker(\overline{\partial}_{b})\subset L^{2}_{(0,q)}(b\Omega)$. $S_{n-1}$ is defined as the orthogonal projection onto $\ker(\overline{\partial}_{b}^{*})^{\perp}$ (equivalently: onto the range of $\overline{\partial}_{b}$). 

The projections $P_{q}$ and $S_{q}$ are said to be continuous in Sobolev--$s$ norms if they map $W^{s}_{(0,q)}(\Omega)$ and $W^{s}_{(0,q)}(b\Omega)$, respectively, to themselves, for $s\geq 0$. By the closed graph theorem, the maps are then automatically continuous (that the graphs are closed follows form continuity in $L^{2}_{(0,q)}(\Omega)$ and $L^{2}_{(0,q)}(b\Omega)$, respectively): there are constants $C_{s}$ and $C_{s}^{\prime}$ such that,
\begin{equation}\label{est0}
 \|P_{q}u\|_{\Omega,s} \leq C_{s}\|u\|_{\Omega,s}\;\;u\in W^{s}_{(0,q)}(\Omega)\;,
\end{equation}
and 
\begin{equation}\label{est0a}
 \|S_{q}u\|_{b\Omega,s} \leq C_{s}^{\prime}\|u\|_{b\Omega,s}\;;\;u\in W^{s}_{(0,q)}(b\Omega)\;,
\end{equation}
respectively. $j_{q}$ is said to be continuous in $W^{s}$--norm if the following estimate holds
\begin{equation}\label{est1}
\|u\|_{\Omega,s} \leq C_{s}(\|\overline{\partial}u\|_{\Omega,s} + \|\overline{\partial}^{*}u\|_{\Omega,s})\;;\;u\in dom(\overline{\partial})\cap dom(\overline{\partial}^{*})\;;\;1\leq q\leq (n-1)\;, 
\end{equation}
with the obvious modifications for $q=0$ and $q=n$. Likewise, $j_{b,q}$ is said to be continuous in $W^{s}$--norm on the boundary if we have the estimate
\begin{equation}\label{est2}
\|u\|_{b\Omega,s} \leq C_{s}(\|\overline{\partial}_{b}u\|_{b\Omega,s} + \|\overline{\partial}_{b}^{*}u\|_{b\Omega,s})\;;\;u\in dom(\overline{\partial}_{b})\cap dom(\overline{\partial}_{b}^{*})\;;\;1\leq q\leq (n-2)\;,
\end{equation}
again with the obvious modifications for $q=0$ and $q=(n-1)$. Both estimates hold in $L^{2}$, i.e. for $s=0$. Note that it is part of the definition that estimates \eqref{est1} and \eqref{est2}  are genuine, i.e. if the right hand side is finite, then so is the left hand side; no \emph{a priori} smoothness of $u$ is assumed.

The embedding $j_{q}$ is compact if and only if the $\overline{\partial}$--Neumann operator $N_{q}$ is compact (\cite{Straube10}, Proposition 4.2). Similarly, $j_{q}$ is subelliptic if and only if $N_{q}$ is subelliptic (\cite{Kohn79}, Theorem 1.13; proof of Theorem 3.4 in \cite{Straube10}). These properties then also propagate up the Sobolev scale. Related statements hold for the operators on the boundary (\cite{BiardStraube17}, Lemma 5; \cite{Kohn85}, (2.1)--(2.2); \cite{KohnNirenberg65}). The purpose of this note is to investigate what happens when only a weaker regularity property of $j_{q}$ or $j_{b,q}$ is assumed, namely regularity in Sobolev spaces. The answer is given by the equivalences in Theorem \ref{equivalent} below. These equivalences complement those between Sobolev continuity of the projection operators (on three consecutive levels) and Sobolev continuity of the $\overline{\partial}$--Neumann operators and complex Greeen operators, respectively (at the middle level), shown in \cite{BoasStraube90} and \cite{HPR15}, as well as those mentioned in Remark 1.2 below.

\begin{thm}\label{equivalent}
 Let $\Omega$ be a smooth bounded pseudoconvex domain in $\mathbb{C}^{n}$.
 
 1) Let $n\geq 2$, $0\leq q\leq n$, $s\geq 0$. Then $j_{q}$ is continuous in $\|\cdot\|_{\Omega,s}$ if and only $P_{q}$ is.
 
 2) Let $n\geq 3$, $0\leq q\leq (n-1)$, $s\geq 0$. Then $j_{b,q}$ is continuous in $\|\cdot\|_{b\Omega,s}$ if and only if $S_{q}$ is.
 
 \noindent Both statements remain true with continuity in Sobolev norms replaced by continuity in $C^{\infty}_{(0,q)}(\overline{\Omega})$ and $C^{\infty}_{(0,q)}(b\Omega)$, respectively.
 \end{thm}

We require $n\geq 3$ in part 2) of Theorem \ref{equivalent} because one of the main tools in the proof, the weighted theory for $\overline{\partial}_{b}$, is incomplete when $n=2$. For example, the arguments in \cite{KohnNicoara06} only allow the authors to obtain regularity of the weighted Szeg\"{o} projection with loss of a derivative (see \cite{KohnNicoara06}, p. 264).

It is worth noting that $P_{n-1}$ is always regular, in view of the formula $P_{n-1}=I-\overline{\partial}^{*}N_{n}\overline{\partial}$ and the elliptic gain of two derivatives for $N_{n}$. Theorem \ref{equivalent} then implies that $j_{n-1}$ is always regular. In particular, $j_{1}$ is always regular $\mathbb{C}^{2}$, notwithstanding the fact that $N_{1}$ need not be (\cite{Barrett92, Christ96}).

\underline{\emph{Remark 1.2:}} It is implicit in \cite{BoasStraube90}, that for $0\leq q\leq (n-1)$, the Bergman projection $P_{q}$ is regular (in Sobolev norms or in $C^{\infty}(\overline{\Omega})$) if and only if the canonical solution operator $\overline{\partial}^{*}N_{q+1}$ for $\overline{\partial}$ is regular on $ker(\overline{\partial})$. Indeed, formula (2) in \cite{BoasStraube90}, with $q$ replaced by $(q+1)$, immediately implies that regularity of $\overline{\partial}^{*}N_{q+1}$ gives regularity of $P_{q}$. The other direction is standard: if $P_{q}$ is regular, then the canonical solution to $\overline{\partial}$, given by a (regular) weighted solution followed by the projection onto $ker(\overline{\partial})^{\perp}$, is also regular. Likewise, it is implicit in \cite{HPR15} that for $0\leq q\leq (n-2)$, the Szeg\"{o} projection $S_{q}$ is regular if and only if the canonical solution operator $\overline{\partial}_{b}^{*}G_{q+1}$ for $\overline{\partial}_{b}$ is regular. The argument is analogous to that for the Bergman projections, but using formula (6) in \cite{HPR15} (see \eqref{eq2} below).

\section{Proof of Theorem \ref{equivalent}}\label{proof}

 \begin{proof}[Proof of Theorem \ref{equivalent}]
The proofs of 1) and 2) are completely analogous, using formulas and ideas from \cite{BoasStraube90} and their analogues form \cite{HPR15}, respectively. We will go through the details of 2), as the boundary operators and the associated weighted theory are perhaps less familiar.

Assume that $S_{q}$ is regular in $s$--norm, for some $s\geq 0$. We must show \eqref{est2}. Consider the case $1\leq q\leq (n-2)$ first. The Hodge decomposition (\cite{ChenShaw01}, Theorem 9.4.2) gives for $u\in dom(\overline{\partial}_{b})\cap dom(\overline{\partial}_{b}^{*})$:
\begin{equation}\label{Hodge}
u = \overline{\partial}_{b}^{*}\overline{\partial}_{b}G_{q}u + \overline{\partial}_{b}\overline{\partial}_{b}^{*}G_{q}u = \overline{\partial}_{b}^{*}G_{q+1}(\overline{\partial}_{b}u) + \overline{\partial}_{b}G_{q-1}(\overline{\partial}_{b}^{*}u) 
= \overline{\partial}_{b}^{*}G_{q+1}(\overline{\partial}_{b}u) + G_{q}\overline{\partial}_{b}(\overline{\partial}_{b}^{*}u)\;.
\end{equation}
So we `only' need bounds in $s$--norm for the operators $\overline{\partial}_{b}^{*}G_{q+1}$ and $G_{q}\overline{\partial}_{b}$ (which are continuous in $L^{2}$). The argument is well known in the case of the $\overline{\partial}$--Neumann operator and relies on Kohn's weighted theory (see for example \cite{Straube10}, Section 5.3). In the case at hand, one replaces this theory by the appropriate weighted theory on the boundary due to Nicoara (\cite{Nicoara06}) and Harrington--Raich (\cite{HR11}; see also \cite{HPR15}). The subscript $t$  denotes operators in the weighted norm; properties of these operators that we need can be found in \cite{HR11}, Theorem 1.2.

The form $v=\overline{\partial}_{b}^{*}G_{q+1}(\overline{\partial}_{b}u)$ is the canonical solution (the solution orthogonal to $\ker(\overline{\partial}_{b})$) to the equation $\overline{\partial}_{b}v=\overline{\partial}_{b}u$. This solution is also given by the weighted canonical solution followed by the orthogonal projection onto $\ker(\overline{\partial}_{b})^{\perp}$, that is by $(I - S_{q})\overline{\partial}_{b,t}^{*}G_{t,q+1}(\overline{\partial}_{b}u)$. When $q\leq (n-3)$, given $s\geq 0$, we can take $t$ large enough so that $\overline{\partial}_{b,t}^{*}G_{t,q+1}$ is continuos in $s$--norm. Thus we obtain
\begin{equation}\label{est3}
 \|\overline{\partial}_{b}^{*}G_{q+1}(\overline{\partial}_{b}u)\|_{b\Omega,s} = \|(I - S_{q})\overline{\partial}_{b,t}^{*}G_{t,q+1}(\overline{\partial}_{b}u)\|_{b\Omega,s} \leq C_{s}\|\overline{\partial}_{b}u\|_{b\Omega,s} \;.
\end{equation}
The second inequality comes form the assumption on $S_{q}$. When $q=(n-2)$, we use $G_{t,q}\overline{\partial}_{b,t}^{*}$ as the weighted canonical solution operator, and again \cite{HR11}, Theorem 1.2, part (vi), to obtain the bound in \eqref{est3}. That $G_{t,q}\overline{\partial}_{b}^{*}$ is the canonical solution operator follows form the general (weighted) $L^{2}$--theory, for example form the weighted analogue of formula (3) in \cite{BoasStraube91b}.

We now estimate $G_{q}\overline{\partial}_{b}(\overline{\partial}_{b}^{*}u)$. The idea comes from \cite{BoasStraube90}, page 29. First note that $\overline{\partial}_{b}^{*}G_{q}=(I-S_{q-1})(\overline{\partial}_{b}^{*}G_{q})S_{q}$ (both sides agree on $\ker(\overline{\partial}_{b})$ and annihilate $\ker(\overline{\partial}_{b})^{\perp}$). Again writing the canonical solution operator in terms of its weighted analogue and the projection onto $\ker(\overline{\partial}_{b})^{\perp}$, and then taking the (unweighted) adjoint, we obtain
\begin{multline}\label{eq4}
G_{q}\overline{\partial}_{b} = (\overline{\partial}_{b}^{*}G_{q})^{*} = S_{q}(\overline{\partial}_{b,t}^{*}G_{t,q})^{*}(I-S_{q}) \\
= S_{q}F_{t}(\overline{\partial}_{b,t}^{*}G_{t,q})_{t}^{*}F_{t}^{-1}(I-S_{q-1})=
S_{q}F_{t}(G_{t,q}\overline{\partial}_{b})F_{t}^{-1}(I-S_{q-1})\;,
\end{multline}
where both $F_{t}$ and its inverse $F_{t}^{-1}$ are self-adjoint pseudodifferential operators of order zero; see Section 1.3 in \cite{HPR15} or \cite{Nicoara06}, Corollary 4.6, for the microlocal setup of the weighted norms. Again choosing $t$ large enough (so that $G_{t,q}\overline{\partial}_{b,t}$ is continuous in $\|\cdot\|_{b\Omega,s}$) gives
\begin{equation}\label{est5}
 \|G_{q}\overline{\partial}_{b}(\overline{\partial}_{b}^{*}u)\|_{b\Omega,s}=\|S_{q}F_{t}(G_{t,q}\overline{\partial}_{b})F_{t}^{-1}\underbrace{(I-S_{q-1})(\overline{\partial}_{b}^{*}u)}_{=\overline{\partial}_{b}^{*}u}\|_{b\Omega,s} \leq C_{s}\|\overline{\partial}_{b}^{*}u\|_{b\Omega,s}\| \;.
\end{equation}
This concludes the case $1\leq q\leq (n-2)$.

When $q=0$ and $u\in\ker(\overline{\partial})^{\perp}$, we have $\|u\|_{b\Omega,s}=\|\overline{\partial}_{b}^{*}G_{1}(\overline{\partial}_{b}u)\|_{b\Omega,s}$; now observe that the argument in \eqref{est3} also works for $q=0$. When $q=(n-1)$, write $u=\overline{\partial}_{b}\overline{\partial}_{b}^{*}G_{n-1}u=\overline{\partial}_{b}G_{n-2}(\overline{\partial}_{b}^{*}u)= (G_{n-1}\overline{\partial}_{b})(\overline{\partial}_{b}^{*}u)$. The case $q=(n-1)$ for $\overline{\partial}_{b}^{*}G_{n-1} = G_{n-2}\overline{\partial}_{b}^{*}$ does not seem to be formally covered in \cite{ChenShaw01}, Theorem 9.4.2, but it follows for example from formula (3) in \cite{BoasStraube91b} (which does cover $q=(n-1)$). That $\overline{\partial}_{b}G_{n-2} = G_{n-1}\overline{\partial}_{b}$ is covered in \cite{ChenShaw01}. Now we use \eqref{eq4} and on the expression on the right hand, replace the weighted canonical solution operator $G_{t,n-1}\overline{\partial}_{b}$ to with $\overline{\partial}_{b}G_{t,n-2}$, and use \cite{HR11}, Theorem 1.2, part (vii) (the commutator relations between the Green operators and $\overline{\partial}_{b}$ in \cite{ChenShaw01}, Theorem 9.4.2, also hold for the weighted operators; alternatively, use again the weighted version of formula (3) ion \cite{BoasStraube91b}).

Now we assume that $j_{b,q}$ is regular in $\|\cdot\|_{b\Omega,s}$ for some $s\geq 0$. We estimate $\|S_{q}u\|_{b\Omega,s}$ by $\|u\|_{b\Omega,s}$ for $u\in C^{\infty}_{(0,q)}(b\Omega)$; by density, the estimate then carries over to $u\in W^{s}_{(0,q)}(b\Omega)$. First, let $0\leq q\leq (n-2)$. Formula (6) in \cite{HPR15}, with $(q-1)$ replaced by $q$, gives
\begin{equation}\label{eq2}
S_{q}u = F_{t}S_{q,t}F_{t}^{-1}u - \overline{\partial}_{b}^{*}G_{q+1}[\overline{\partial}_{b},F_{t}]S_{q,t}F_{t}^{-1}u\;;
\end{equation}
here $F_{t}$ is as above. Technically, for (6) in \cite{HPR15} to apply, one needs $q\leq (n-3)$, but it is easy to see that their calculation also works for $q=(n-2)$ (i.e. `their' $q$ equals $(n-1)$). If we take $t$ big enough, the $s$--norm of the first term is dominated by $\|u\|_{b\Omega,s}$ as desired. For the second term, we use the assumption on $j_{q}$:
\begin{multline}\label{est6}
\|\overline{\partial}_{b}^{*}G_{q+1}[\overline{\partial}_{b},F_{t}]S_{q,t}F_{t}^{-1}u\|_{b\Omega,s} \leq C_{s}\|\overline{\partial}_{b}\overline{\partial}_{b}^{*}G_{q+1}[\overline{\partial}_{b},F_{t}]S_{q,t}F_{t}^{-1}u\|_{b\Omega,s} \\
= C_{s}\|[\overline{\partial}_{b},F_{t}]S_{q,t}F_{t}^{-1}u\|_{b\Omega,s}\;.
\end{multline}
In the equality, we have used that $[\overline{\partial}_{b},F_{t}]S_{q,t}F_{t}^{-1}u = \overline{\partial}_{b}F_{t}S_{q,t}F_{t}^{-1}u$ is in the range of $\overline{\partial}_{b}$, and that on such forms, $\overline{\partial}_{b}\overline{\partial}_{b}^{*}G_{q+1}$ acts as the identity. Taking into account that the commutator $[\overline{\partial}_{b},F_{t}]$ is of order zero shows that the right hand side of \eqref{est6} is also bounded by $\|u\|_{b\Omega,s}$.

Now let $q=(n-1)$. $S_{n-1}$ agrees with the projection $S^{\prime}_{n-1}$ in \cite{HPR15}. The first line in formula (7) from \cite{HPR15}, with $q$ replaced by $(n-2)$, then says $S_{n-1}u = S_{n-1,t} + \overline{\partial}_{b}G_{n-2}\overline{\partial}_{b}^{*}(I- S_{n-1,t})u$. $\overline{\partial}_{b}G_{n-2}\overline{\partial}_{b}^{*}(I- S_{n-1,t})u$ is in the range of $\overline{\partial}_{b}$, and we can apply the assumption on $j_{b,n-1}$ to it. Doing so gives
\begin{multline}\label{eq7}
\|S_{n-1}u\|_{b\Omega,s} \leq C_{s}(\|S_{n-1,t}u\|_{b\Omega,s}+\|\overline{\partial}_{b}^{*}\overline{\partial}_{b}G_{n-2}\overline{\partial}_{b}^{*}(I-S_{n-1,t})u\|_{b\Omega,s}) \\
= C_{s}(\|S_{n-1,t}u\|_{b\Omega,s}+\|\overline{\partial}_{b}^{*}(I-S_{n-1,t})u\|_{b\Omega,s})\;\;\;\;\; \\
= C_{s}(\|S_{n-1,t}u\|_{b\Omega,s}+\|(\overline{\partial}_{b}^{*}-\overline{\partial}_{b,t}^{*})(I-S_{n-1,t})u\|_{b\Omega,s}) \lesssim\|u\|_{b\Omega,s}\;,
\end{multline}
for $t$ big enough. We have used here that $\overline{\partial}_{b,t}^{*}$ annihilates $(I-S_{n-1,t})u$, and that the difference $(\overline{\partial}_{b}^{*}-\overline{\partial}_{b,t}^{*}) =\overline{\partial}_{b}^{*}F_{t}F_{t}^{-1}-F_{t}\overline{\partial}_{b}^{*}F_{t}^{-1}=[\overline{\partial}_{b}^{*},F_{t}]F_{t}^{-1}$ is of order zero (compare \cite{HPR15, BoasStraube90}). That $S_{n-1,t}=\overline{\partial}_{b}G_{n-2,t}\overline{\partial}_{b,t}^{*}$ is continuous in $\|\cdot\|_{b\Omega,s}$ for $t$ big enough is not stated explicitly in the theorems in \cite{HR11,HPR15}, but is in line 6 from the bottom of page 362 in \cite{HPR15}, in light of the estimate for $G_{n-2,t}\overline{\partial}_{b,t}^{*}$ in Proposition 3.1 (\cite{Raich24}; we allow the constant in the estimate to depend on $t$).

This completes the proof of (2) in Theorem \ref{equivalent}. The corresponding statement with $C^{\infty}_{(0,q)}(\overline{\Omega})$ and $C^{\infty}_{(0,q)}(b\Omega)$ follows along the same lines, using that the Sobolev norms $\{\|\cdot\|_{k}\}_{k\in\mathbb{N}}$ and $\{\|\cdot\|_{b\Omega,k}\}_{k\in\mathbb{N}}$ induce the Fr\'{e}chet topology in the  respective $C^{\infty}$ spaces.

Finally, as mentioned above, part (1) and its analogue in the $C^{\infty}$--spaces follows almost verbatim, but using the formulas from \cite{BoasStraube90}. In this case, the (pseudodifferential) operator $F_{t}$ is simply multiplication by the weight function $e^{-t|z|^{2}}$. We leave the details to the reader. This completes the proof of Theorem \ref{equivalent}.
\end{proof}

\bigskip
\bigskip

\emph{Acknowledgement:} I am indebted to Andy Raich for helpful correspondence on his paper \cite{HPR15}.

\vskip .5cm

\providecommand{\bysame}{\leavevmode\hbox to3em{\hrulefill}\thinspace}


\begin{thebibliography}{10}

\bibitem{Barrett92}
Barrett, David E., Behavior of the Bergman projection on the
Diederich-Forn\ae ss worm, \emph{Acta Math.} \textbf{168} (1992), 1--10.

%\bibitem{BarrettSahutoglu12}
%Barrett, David.~E. and \c{S}ahuto\u{o}lu, Sönmez,
%Irregularity of the Bergman projection on worm domains in $\mathbb{C}^{n}$.
%\emph{Michigan Math. J.} \textbf{61}, no.1 (2012), 187–198.

\bibitem{BiardStraube17}
Biard, S\'{e}verine and Straube, Emil J., $L^{2}$--Sobolev theory for the complex Green operator, \emph{Internat. J. Math.} {\bf 28}, no. 9 (2017), 1740006, 31 pp.

%\bibitem{BiardStraube19}
%\bysame, Estimates for the complex Green operator: symmetry, percolation, and interpolation, \emph{Trans. Amer. Math. Soc.} {\bf 371}, no. 3 (2019), 2003–2020.

\bibitem{BoasStraube90}
Boas, Harold P. and Straube, Emil J., Equivalence of regularity for the Bergman projection and the $\overline{\partial}$-Neumann operator, \emph{Manuscripta Math.} \textbf{67} (1990), 25--33.

\bibitem{BoasStraube91b}
\bysame, Sobolev estimates for the complex Green operator on a class of weakly pseudoconvex boundaries, \emph{Comm. Partial Differential Equations} {\bf 16}, no. 10 (1991), 1573–1582. 

%\bibitem{Catlin87b}
%Catlin, David W., Subelliptic estimates for the $\overline{\partial}$-Neumann problem, \emph{Ann. Math. (2)}, \textbf{126} (1987), 131--191.

\bibitem{Boggess91}
Boggess, Albert, \emph{CR-Manifolds and the Tangential Cauchy-Riemann Complex}, Studies in Advanced Mathematics, CRC Press 1991.

\bibitem{ChenShaw01}
Chen, So-Chin and Shaw, Mei-Chi, \emph{Partial Differential Equations in
Several Complex Variables}, Studies in Advanced Mathematics 19, Amer. Math. Soc./International Press, 2001.

\bibitem{Christ96}
Christ, Michael, Global $C^{\infty}$ irregularity of the
$\overline{\partial}$-Neumann problem for worm domains, \emph{J. Amer.
Math. Soc.} \textbf{9}, Nr. 4 (1996), 1171--1185.

%\bibitem{Diaz86}
%Diaz, Ricardo L., Necessary conditions for subellipticity of $\Box_{b}$ on pseudoconvex boundaries. \emph{Comm. Partial Differential Equations} \textbf{11} (1) (1986) 1–61.

\bibitem{HPR15}
Harrington, Phillip S., Peloso, Marco M., and Raich, Andrew S., Regularity equivalence of the Szegö projection and the complex Green operator, \emph{Proc. Amer. Math. Soc.} \textbf{143} no.1 (2015), 353–367.

\bibitem{HR11}
Harrington, Phillip S. and Raich, Andrew, Regularity results for $\overline{\partial}_{b}$ on CR-manifolds of hypersurface type, \emph{Comm. Partial Differential Equations} \textbf{36}, no.1 (2011), 134–161.

%\bibitem{Khanh16}
%Khanh, Tran Vu, Equivalence of estimates on a domain and its boundary, \emph{Vietnam J. Math.} \textbf{44}, no. 1 (2016), 29–48.

%\bibitem{KhanhZampieri13}
%\bysame  and Zampieri, Giuseppe, Regularity at the boundary and tangential regularity of solutions of the Cauchy-Riemann system, \emph{Pacific J. Math.} \textbf{265} no.2 (2013), 491–498.

%\bibitem{Koenig04}
% Koenig, Kenneth D., A parametrix for the $\overline{\partial}$--Neumann problem on pseudoconvex domains of finite type, \emph{J. Funct. Anal.} \textbf{216}, no.2 (2004), 243–302.

\bibitem{Kohn79}
Kohn, J.~J., Subellipticity of the $\overline{\partial}$-Neumann problem on pseudoconvex domains: sufficient conditions, \emph{Acta Math.} \textbf{142} (1979), 79--122.

%\bibitem{Kohn79b}
%Kohn, J.~J., Boundary regularity of $\overline{\partial}$, \emph{Recent developments in several complex variables} (Proc. Conf., Princeton Univ., Princeton, N. J., 1979), pp. 243–260, Ann. of Math. Stud., No. 100, Princeton Univ. Press, Princeton, NJ, 1981.

\bibitem{Kohn85}
Kohn, J. J., Estimates for $\overline{\partial}_{b}$ on pseudoconvex CR manifolds, \emph{Pseudodifferential operators and applications} (Notre dame, 1984), Proc. Sympos. Pure Math., 43, Amer. Math. Soc., Providence, RI, 1985, 207--217.

%\bibitem{Kohn02}
%\bysame, Superlogarithmic estimates on pseudoconvex domains and CR manifolds,  \emph{Ann.of Math.}(2) \textbf{156} (2002), 213--248.

\bibitem{KohnNicoara06}
Kohn, J.~J. and Nicoara, A.~C., The $\overline\partial\sb b$ equation on weakly pseudo-convex CR manifolds of dimension 3, \emph{J. Funct. Anal.} \textbf{230}, no.2  (2006), 251--272.

\bibitem{KohnNirenberg65}
Kohn, J.~J. and Nirenberg, L., \emph{Non-coercive boundary value
problems}, Comm. Pure Appl. Math. \textbf{18} (1965), 443--492.

\bibitem{Nicoara06}
Nicoara, Andreea C., Global regularity for $\overline\partial\sb b$ on weakly pseudoconvex CR manifolds, \emph{Adv. Math.} \textbf{199}, no.2  (2006), 356--447.

\bibitem{Raich24}
Raich, Andrew S., private communication.

%\bibitem{RaichStraube08}
%Raich, Andrew S. and Straube, Emil J., Compactness of the complex Green operator, \emph{Math. Res. Lett.} \textbf{15}, no.~4 (2008), 761--778.

\bibitem{Straube10}
Straube, Emil J., \emph{Lectures on the $L^2$-Sobolev theory of the $\overline{\partial}$--Neumann problem}, ESI Lectures in Mathematics and Physics, European Mathematical Society (EMS), Z\"{u}rich 2010.



\end{thebibliography}
\end{document}